\documentclass[graybox]{svmult}

\usepackage{type1cm}        
\usepackage{makeidx}         
\usepackage{graphicx}        
\usepackage{multicol}        
\usepackage[bottom]{footmisc}

\usepackage{newtxtext}       %
\usepackage[varvw]{newtxmath}       

\makeindex             

\begin{document}

\title*{On the Lebesgue structure of the distribution of a random variable defined by continued \(A_2$-fractions}
\titlerunning{On the Lebesgue structure of the distribution of a random variable} 
\author{Pratsiovytyi Mykola\orcidID{0000-0001-6130-9413}, \\ Makarchuk Oleg \orcidID{0000-0002-1001-8568} and \\Karvatskyi Dmytro\orcidID{0000-0002-6873-6271}}
\institute{Pratsiovytyi Mykola \at Dragomanov Ukrainian State University, Pyrohova Str. 9, Kyiv, 01601, Ukraine
\at Institute of Mathematics of NAS of Ukraine, Tereschenkivska Str. 3, Kyiv, 01601, Ukraine \newline \email{pratsiovytyi@imath.kiev.ua}
\and Makarchuk Oleg \at Institute of Mathematics of NAS of Ukraine, Tereschenkivska Str. 3, Kyiv, 01601,
Ukraine \newline \email{makolpet@gmail.com}
\and Karvatskyi Dmytro \at Institute of Mathematics of NAS of Ukraine, Tereschenkivska Str. 3, Kyiv, 01601,
Ukraine \newline \email{karvatsky@imath.kiev.ua}}
%
%
\maketitle

\abstract*{In this paper, we study the Lebesgue structure of the distribution of a random variable given in terms of a continued fraction with a two-symbol alphabet $\{\frac{1}{2}; 1\}$, also known as $A_2$-fractions. We establish necessary and sufficient conditions for the distribution to be discrete and provide some sufficient conditions for its singularity. We also explore non-trivial metric properties of the $A_2$-representation with the specified alphabet.}

\abstract{In this paper, we study the Lebesgue structure of the distribution of a random variable given in terms of a continued fraction with a two-symbol alphabet $\{\frac{1}{2}; 1\}$, also known as $A_2$-fractions. We establish necessary and sufficient conditions for the distribution to be discrete and provide some sufficient conditions for its singularity. We also explore non-trivial metric properties of the $A_2$-representation with the specified alphabet.
}

\section{Introduction}

For a positive real number sequence \((a_{n})$ and a natural number \(k$, we define a finite continued fraction as follows:
$$[a_{1};a_{2};\dots;a_{k}] = \cfrac{1}{a_1
          + \cfrac{1}{a_2
          + \dots + \cfrac{1}{a_{k-1} + \cfrac{1}{a_{k}}}}}. $$

\newpage

Let
\begin{equation}\label{frac}
[a_{1};a_{2};\dots;a_{k};\dots] = \lim_{k \to +\infty}[a_{1};a_{2};\dots;a_{k}]
\end{equation}
provided that the above limit exists. It was shown independently by Seidel and Stern \cite{seid,stern} that the limit in \eqref{frac} exists if and only if the series $\sum_{n=1}^{+\infty}a_{n}$ is divergent.

Let $(\chi_{n})$ be a sequence of independent and identically distributed random variables that take $k$ distinct natural values with probabilities $p_{kn}$, respectively. The Lebesgue structure of the distribution of a random variable $\chi = [\chi_{1};\chi_{2};...;\chi_{n};...]$ was studied in \cite{prlan}, where it was proven that $\chi$ is discrete if and only if
$$\prod_{n=1}^{+\infty}\max_{k \in \mathbb{N}}\{p_{kn}\} > 0,$$
and otherwise, $\chi$ is singular. By other methods, an analogous result was obtained in \cite{prtorlan}, where the properties of the distribution spectrum of $\chi$ were also investigated. Various classes of singular distributions connected with continued fractions were constructed in \cite{letac, lyons, sim}.
 
It is well known from \cite{prats} that for any $t \in [\frac{1}{2}, 1]$, there exists a sequence $(b_{n})$ with $b_n \in \{ \frac{1}{2}, 1\}$ for all $n \in \mathbb{N}$, such that
$$t = [b_{1}; b_{2}; \dots; b_{n}; \dots].$$
This expression is called the $A_2$-representation with the alphabet $\{\frac{1}{2}, 1\}$. A countable subset of $[\frac{1}{2}, 1]$ has two distinct $A_2$-representations of the form
$$[b_{1}; \dots; b_{n}; \frac{1}{2};(\frac{1}{2};1)] = [b_{1}; \dots; b_{n}; 1; (1; \frac{1}{2})],$$
where parentheses denote the period of a given continued fraction. Numbers possessing the above property are called $A_2$-binary. In contrast, numbers in the interval $[\frac{1}{2}, 1]$ that are not $A_2$-binary and have a unique $A_2$-representation are termed $A_2$-unary. A detailed analysis of the topological and metric properties of the $A_2$-representation is provided in \cite{prats}.  

Let $(\xi_{n})$ be a sequence of independent discretely distributed random variables that take values $\frac{1}{2}$ or $1$ with probabilities $p_{(\frac{1}{2})n}$ and $p_{(1)n}$, respectively.  
The Lebesgue structure of the distribution of the random variable $\xi = [\xi_{1}; \xi_{2}, \dots, \xi_{n}, \dots]$ was studied in \cite{kur} in the case $p_{(\frac{1}{2})n} = \gamma$ for each natural $n$ and some $\gamma \in [0, 1]$.  
Using nontrivial methods, it was shown in \cite{kur, prkur} that the distribution of $\xi$ is degenerate when $\gamma \in \{0, 1\}$ and singular when $\gamma \in (0, 1)$.  
In \cite{prkur}, certain fractal properties of the distribution of $\xi$ for the corresponding case were also considered.

The present paper is devoted to the study of the Lebesgue structure of the random variable $\xi$.  
Necessary and sufficient conditions for the discreteness, as well as certain sufficient conditions for the singularity of the distribution of the random variable $\xi$, have been established. The corresponding results extend and deepen the findings of \cite{kur, prkur}.

\section{On a metric property of $A_{2}$-representation}
\label{sec:2}

Let $ W(\frac{1}{2}; 1) $ denote the set of all sequences $ (t_{n}) $ such that $ t_{k} \in \{\frac{1}{2}, 1\} $ for every natural number $k$.  
For a given sequence $ (a_{n}) \in W(\frac{1}{2}; 1) $ and the corresponding number $x = [a_{1}; a_{2}; \dots; a_{k}; \dots]$ we define two sequences, $q_{n}(x)$ and $p_{n}(x)$, by
$$
q_{n}(x) = q_{n}(a_{1}; a_{2}; \dots; a_{n}) = a_{n} q_{n-1}(a_{1}; a_{2}; \dots; a_{n-1}) + q_{n-2}(a_{1}; a_{2}; \dots; a_{n-2}) ,
$$
$$
p_{n}(x) = p_{n}(a_{1}; a_{2}; \dots; a_{n}) = a_{n} p_{n-1}(a_{1}; a_{2}; \dots; a_{n-1}) + p_{n-2}(a_{1}; a_{2}; \dots; a_{n-2}) ,
$$
with initial conditions $p_{-1} = 1, q_{-1} = 0, p_{0} = 0, q_{0} = 1$.

It is well known from \cite{chin, prkur} that
\begin{equation}\label{eq1}
\lim_{n \to +\infty} \frac{p_{n}(x)}{q_{n}(x)} = x, \quad p_{n}(x) q_{n-1}(x) - p_{n-1}(x) q_{n}(x) = (-1)^{n+1} \quad \forall n \in \mathbb{N},
\end{equation}

\begin{equation}\label{eq2}
\frac{q_{n}(x)}{q_{n-1}(x)} = a_{n-1} + [a_{n-2}; \dots; a_{1}] \in [1; 2] \quad \forall n \in \mathbb{N},
\end{equation}

\begin{equation}\label{eq3}
q_{n}(x) \le \frac{2\sqrt{17}}{17} \left( \left( \frac{1+\sqrt{17}}{4} \right)^{n+1} - \left( \frac{1-\sqrt{17}}{4} \right)^{n+1} \right) = w_{n} \quad \forall n \in \mathbb{N}.
\end{equation}

For convenience, we will use simplified notation $ q_{n} = q_{n}(x) = q_{n}(a_{1}; a_{2}; \dots; a_{n}) $ and $ p_{n} = p_{n}(x) = p_{n}(a_{1}; a_{2}; \dots; a_{n}) $) where it is necessary. We also write
$$
[x, y]^{*} = [\min(x, y), \max(x, y)], \quad x, y \in \mathbb{R}.
$$

It was proven in \cite{prats} that for an arbitrary finite sequence $ (\gamma_{n}) \in W(\frac{1}{2}; 1) $ and the corresponding number $ x = [\gamma_{1}; \gamma_{2}; \dots; \gamma_{n}] $, the set
$$
\{ [\gamma_1; \gamma_2; \dots; \gamma_n; \sigma_1; \sigma_2; \dots; \sigma_k; \dots ] \mid \sigma_k \in \{\frac{1}{2}, 1\} \; \forall k \in \mathbb{N} \}
$$
is a closed interval of the form
\begin{equation}
\label{eq4}
\Delta_{n}(x) = \Delta_{n}(\gamma_1; \gamma_2; \dots; \gamma_n) = \left[ \frac{p_n + \frac{1}{2} p_{n-1}}{q_n + \frac{1}{2} p_{n-1}}, \frac{p_n + p_{n-1}}{q_n + p_{n-1}} \right]^{*},
\end{equation}
which is called the cylinder of the $n$-th level. Moreover, two distinct cylinders of the same level have at most one common point (the touching point). For $n$-th level cylinders, we will use the symbol $\Delta_{n}$.

For the $A_2$-representation on the measurable space $([0, 1]; \mathcal{B}[0, 1])$, where $\mathcal{B}[a, b]$ is the sigma-algebra of Borel sets on the interval $[a, b]$, we define the left shift operator $T$ by

\begin{equation}
\label{}
T([a_1; a_2; \dots; a_n; \dots]) = [a_2; a_3; \dots; a_{n+1}; \dots].
\end{equation}

From now on, we agree not to use $A_2$-representations with period $(\frac{1}{2}; 1)$ for $A_2$-binary numbers.  
Let $\lambda(\cdot)$ stand for the standard Lebesgue measure, and let $\eta(\cdot)$ denote the Lebesgue–Stieltjes measure corresponding to the absolutely continuous distribution function
$$
h(x) = \frac{\ln(x+1) - \ln(x+2) + \ln\left(\frac{5}{3}\right)}{\ln\left(\frac{10}{9}\right)} \quad \forall x \in [\frac{1}{2}, 1],
$$
with the density
$$
p(x) = \ln\left(\frac{10}{9}\right) (x+1)(x+2) \quad \forall x \in [\frac{1}{2}, 1].
$$

\begin{theorem}
\label{teo1}
    The dynamical system $\left([\frac{1}{2}, 1], \mathcal{B}[\frac{1}{2}, 1], \eta(\cdot), T\right)$ is ergodic.
\end{theorem}

\begin{proof}
    We will show that $T$ preserves the measure $\eta(\cdot)$. To do this, it is sufficient to show that for any interval $[\alpha, \beta] \subseteq [\frac{1}{2}, 1]$, the condition
    $$\eta(T^{-1}([\alpha, \beta])) = \eta(T([\alpha, \beta]))$$
    is fulfilled. Indeed, the relation
    $$T^{-1}([\alpha, \beta]) = \left[\frac{1}{\beta + \frac{1}{2}}, \frac{1}{\alpha + \frac{1}{2}}\right] \cup \left[\frac{1}{\beta + 1}, \frac{1}{\alpha + 1}\right] 
    \Rightarrow$$
    $$\Rightarrow \eta(T^{-1}([\alpha, \beta])) = \eta(T([\alpha, \beta])) \Leftrightarrow$$
    $$\Leftrightarrow f(\beta) - f(\alpha) = f\left(\frac{1}{\alpha + \frac{1}{2}}\right) - f\left(\frac{1}{\beta + \frac{1}{2}}\right)
        + f\left(\frac{1}{\alpha + 1}\right) - f\left(\frac{1}{\beta + 1}\right) \Leftrightarrow$$
    $$\Leftrightarrow f(\beta) + f\left(\frac{1}{\beta + 1}\right) + f\left(\frac{1}{\beta + \frac{1}{2}}\right) = f(\alpha) + f\left(\frac{1}{\alpha + 1}\right) + f\left(\frac{1}{\alpha + \frac{1}{2}}\right)$$
    is correct since for every $t \in [\frac{1}{2}, 1]$
    $$
    \ln\left(\frac{10}{9}\right)\left(f(t) + f\left(\frac{1}{t + 1}\right) + f\left(\frac{1}{t + \frac{1}{2}}\right)\right) =
    \ln(t + 1) - \ln(t + 2) - \ln(t + 2) - \ln(t + 1)
    $$
    $$-\ln(2t + 3) + \ln(t + 1) + \ln(t + 1.5) - \ln(t + \frac{1}{2}) - \ln(2t + 2) +$$
    $$+\ln(t + \frac{1}{2}) + \ln\left(\frac{5}{3}\right) =
    \ln\left(\frac{5}{3}\right) - 2\ln(2).$$

    Next, we show that the transformation $T$ is ergodic, i.e., for any invariant set $W$ such that $T^{-1}(W) = W$, the condition $\eta(W) \in \{0, 1\}$ is fulfilled.  
    Let $[a, b] \subseteq [\frac{1}{2}, 1]$ and $(c_{n}) \in W(\frac{1}{2}; 1)$. Taking into account \eqref{eq1}, we get
    $$
    \lambda\left(T^{-n}([a, b]) \cap \Delta_n(c_1; c_2; \ldots; c_n)\right) = $$ 
    $$ = \lambda\left(\left[\left[c_1; c_2; \ldots; c_n; \frac{1}{a}\right],
    \left[c_1; c_2; \ldots; c_n; \frac{1}{b}\right]\right]^* \right) =
    $$
    $$
    = \left|\frac{p_{n-1} + \frac{1}{b}p_n}{q_{n-1} + \frac{1}{b}q_n} - \frac{p_{n-1} + \frac{1}{a}p_n}{q_{n-1} + \frac{1}{a}q_n}\right| =
    \frac{(b - a)}{(bq_{n-1} + q_n)(aq_{n-1} + q_n)} =
    $$
    $$
    = \lambda([a, b]) \cdot \lambda(\Delta_n) \cdot \frac{2(\frac{1}{2}q_{n-1} + q_n)(q_{n-1} + q_n)}{(bq_{n-1} + q_n)(aq_{n-1} + q_n)} \;\; \forall n \in \mathbb{N}.
    $$
    Next we consider the function $h(z) = \frac{z + \frac{1}{2}}{z + 1} = 1 - \frac{1}{2z + 2}$, that is increasing on the interval $[1, 2]$. According to \eqref{eq2}, we can observe
    $$
    \frac{(\frac{1}{2}q_{n-1} + q_n)(q_{n-1} + q_n)}{(bq_{n-1} + q_n)(aq_{n-1} + q_n)} \geq \frac{\frac{1}{2}q_{n-1} + q_n}{q_{n-1} + q_n} = h\left(\frac{q_n}{q_{n-1}}\right) \geq h(1) = \frac{3}{4}.
    $$
    Applying this to the terms in the measure expression, we arrive at the inequality
    $$\lambda\left(T^{-n}([a, b]) \cap \Delta_n\right) \geq \frac{3}{4}\lambda([a, b]) \cdot \lambda(\Delta_n).$$

    Let $E$ be an invariant set under the transformation $T$ with $\lambda(E) > 0$. Since $T^{-n}(E) = E$ for any $n$, we get
    $$\lambda\left(T^{-n}(E) \cap \Delta_n\right) = \lambda\left(E \cap \Delta_n\right) \geq \frac{3}{4}\lambda(E)\lambda(\Delta_n).$$
    Taking into account Knopp's lemma \cite{daj}, we obtain $\lambda(E) = 1$. Since the measure $\eta(\cdot)$ is equivalent to the measure $\lambda(\cdot)$, it follows that $\eta(E) = 1$.
\end{proof}

\begin{theorem}\label{teo2}
    Let $G$ be defined as
    \begin{equation}\label{g}
    G = -\int_{\frac{1}{2}}^{1}\frac{\ln(x)}{\ln\left(\frac{10}{9}\right)(x+1)(x+2)}dx \approx 0.3337.
    \end{equation}
    For almost all numbers $ x \in [\frac{1}{2}; 1]$ in the sense of Lebesgue measure, the following equality holds:
    \begin{equation}\label{lev2}
    \lim\limits_{n\rightarrow\infty} \sqrt[n]{\lambda(\Delta_n)}=e^{-2G} \approx \frac{1}{2}134.
    \end{equation}
\end{theorem}

\begin{proof}
    Since $ p_{1}(T^{n-1}(x)) = 1 $ and $ p_n(x) = q_{n-1}(T(x)) $ for every natural number $ n $, we get the identity
    \begin{align*}
        \frac{p_n(x)}{q_n(x)} \cdot \frac{p_{n-1}(T(x))}{q_{n-1}(T(x))} \cdot \ldots \cdot \frac{p_{1}(T^{n-1}(x))}{q_{1}(T^{n-1}(x))} = \\
        = \frac{p_n(x)}{q_n(x)} \cdot \frac{p_{n-1}(T(x)) p_{n-2}(T^2(x)) \ldots p_{1}(T^{n-1}(x))}{p_{n}(x) p_{n-1}(T(x)) \ldots p_{2}(T^{n-2}(x))}
        &= \frac{1}{q_n(x)}.
    \end{align*}
    It follows that
    $$
    -\ln(q_n(x)) = \ln\frac{p_n(x)}{q_n(x)} + \ln\frac{p_{n-1}(T(x))}{q_{n-1}(T(x))} + \ldots + \ln\frac{p_{1}(T^{n-1}(x))}{q_{1}(T^{n-1}(x))}.
    $$

    Let $ z \in \Delta_k $, then for some $ z_{1} \in \Delta_k $:
    $$
    \left|\ln(z) - \ln\left(\frac{p_{k}(z)}{q_{k}(z)}\right)\right| = \frac{1}{z_{1}} \cdot \left|z - \frac{p_{k}(z)}{q_{k}(z)}\right| \leq 2\lambda(\Delta_k).
    $$
    Considering \eqref{eq1}, \eqref{eq3}, and \eqref{eq4}, we obtain
    $$
    \lambda(\Delta_k) = \left|\frac{p_n + \frac{1}{2}p_{n-1}}{q_n + \frac{1}{2}q_{n-1}} - \frac{p_n + p_{n-1}}{q_n + q_{n-1}}\right| 
    = \frac{1}{(q_{k-1} + q_k)(q_{k-1} + 2q_k)} \leq 
    $$
    $$
    \leq \frac{1}{(w_{k-1} + w_k)(w_{k-1} + 2w_k)} = c_{k}, \;\; \forall k \in \mathbb{N}.
    $$

    It is clear that the series $ \sum\limits_{n=0}^{+\infty}c_{n} $ converges and has a positive sum $ L $.
    Thus, for each $ j \in \{0, 1, \ldots, n-1\} $, we have
    $$
    \ln\left(\frac{p_{n-j}(T^j(x))}{q_{n-j}(T^j(x))}\right) = \ln(T^j(x)) + \Delta_j(x),
    $$
    where $ \Delta_j(x) \leq 2c_{n-j} $. By Theorem \ref{teo1} and the Birkhoff-Khinchin theorem, for almost all $ x \in [\frac{1}{2}, 1] $ in the sense of the measure $ \eta(\cdot) $, we have
    $$
    \lim_{n \to +\infty}\frac{1}{n} \sum_{k=0}^{n-1}\ln \left(T^{k}(x)\right) = \int_{\frac{1}{2}}^{1}\frac{\ln(x)}{\ln\left(\frac{10}{9}\right)(x+1)(x+2)}dx = -G.
    $$

    Since for every natural number $n$
    $$
    \left|\sum\limits_{k=0}^{n-1}\Delta_k(x)\right| \leq \sum\limits_{k=0}^{n-1}2c_k \leq 2L,
    $$
    for almost all $ x \in [\frac{1}{2}, 1] $, we conclude that
    $$
    \frac{1}{n}\ln(q_n) = -\frac{1}{n}\sum\limits_{k=0}^{n-1}\ln(T^{k}(x)) - \frac{1}{n}\sum\limits_{k=0}^{n-1}\Delta_k(x) \to G \quad (n \to +\infty).
    $$
    It implies that
    \begin{equation}\label{lev1}
        \lim_{n \to +\infty}\sqrt[n]{q_{n}} = e^{G}.
    \end{equation}

    From the inequality \( \frac{1}{2}q_k \leq q_{k-1} \leq q_k \) for each natural number \( k \), we deduce that the measure \( \lambda(\Delta_k) \) satisfies the following bounds:
$$
\lambda(\Delta_k) = \frac{1}{(q_{k-1} + q_k)(q_{k-1} + 2q_k)} \in \left[\frac{1}{6q_k^2}, \frac{1}{3.75q_k^2}\right], \quad \forall k \in \mathbb{N}.
$$
By using these estimations in combination with equality \eqref{lev1}, we obtain \eqref{lev2}.

\end{proof}

\section{The Lebesgue structure of the distribution of $\xi$}

\begin{theorem}\label{teo3}
    The distribution of $\xi$ is discrete if and only if
    \begin{equation}\label{Form3}
        \prod_{n=1}^{+\infty}\max\left\{p_{(\frac{1}{2})n}, p_{(1)n}\right\} > 0.
    \end{equation}
\end{theorem}

\begin{proof}
    Suppose the distribution of $\xi$ is discrete. Then, for some $a \in [\frac{1}{2}, 1]$, we have $P(\xi = a) > 0.$
    If $a$ is $A_2$-binary, there exist sequences $(\gamma_n), (\beta_n) \in W(\frac{1}{2}; 1)$ such that
    $$
        a = [\gamma_1; \gamma_2; \ldots; \gamma_n; \ldots] = [\beta_1; \beta_2; \ldots; \beta_n; \ldots].
    $$
    Then,
    $$
        \prod_{n=1}^{+\infty} p_{(\gamma_n)n} + \prod_{n=1}^{+\infty} p_{(\beta_n)n} = P(\xi = a) > 0.
    $$
    Thus, at least one of the products $\prod_{n=1}^{+\infty} p_{(\gamma_n)n}$ or $\prod_{n=1}^{+\infty} p_{(\beta_n)n}$ is greater than zero, and condition \eqref{Form3} holds. If $a$ is not $A_2$-binary, a similar reasoning applies.

    Conversely, suppose condition \eqref{Form3} holds, i.e., there exists a sequence $(\gamma_n) \in W(\frac{1}{2}; 1)$ such that
    $$
        \prod_{n=1}^{+\infty} p_{(\gamma_n)n} > 0.
    $$
    Assume that $F_\xi(x)$ can be expressed as
    $$
        F_{\xi}(x) = \alpha_d F_d(x) + \alpha_c F_c(x),
    $$
    where $F_d(x)$ and $F_c(x)$ are the discrete and continuous distribution functions, respectively, and $\alpha_d + \alpha_c = 1$ with $\alpha_d < 1$.
    For a natural number $m$, consider the set $A_m$ of numbers having form
    $$
        [b_1; b_2; \ldots; b_m; \gamma_1; \gamma_2; \ldots; \gamma_l; \ldots],
    $$
    where $b_j \in \{\frac{1}{2}, 1\}$ for each $j \in \{1, 2, \ldots, m\}$. Then,
    $$
        \alpha_d \geq \alpha_d F_d(A_m) = F_{\xi}(A_m) = \prod_{n=m+1}^{+\infty} p_{(\gamma_n)n} \to 1 \quad (m \to +\infty).
    $$
    This leads to a contradiction, as $\alpha_d < 1$. Therefore, the distribution must be discrete.
\end{proof}

\begin{theorem}
    Suppose there exist numbers $M \in \mathbb{Z}_{+}$, $L \in \mathbb{N}$, and stochastic vectors $(q_{(\frac{1}{2})j}, q_{(1)j})$, where $j \in \{1, 2, \ldots, L\}$, such that
    $$
        p_{(\frac{1}{2})(M+jL+i)} = q_{(\frac{1}{2})i} \quad \forall j \in \mathbb{Z}_{+}, \; i \in \{1, 2, \ldots, L\}.
    $$
    If the following conditions hold:
    $$
        \exists j \in \{1, 2, \ldots, L\}: q_{(\frac{1}{2})j} \in (0, 1),
    $$
    \begin{equation}\label{main}
        \prod_{j=1}^{L} (q_{(\frac{1}{2})j})^{\frac{\ln\frac{25}{24}}{\ln\frac{10}{9}}} (q_{(1)j})^{1 - \frac{\ln\frac{25}{24}}{\ln\frac{10}{9}}} \neq e^{-2GL},
    \end{equation}
    where $G$ is defined by equality \eqref{g}, then the distribution of $\xi$ is singular.
\end{theorem}

\begin{proof}
    By Theorem \ref{teo3}, the function $F_{\xi}(x)$ is continuous. Suppose that the Lebesgue decomposition of $F_{\xi}(x)$ contains an absolutely continuous component. Let us distinguish two subcases:

    \textbf{Subcase 1:} $q_{(\frac{1}{2})j} \in (0, 1)$ for each $j \in \{1, 2, \ldots, L\}$. Now, $\tau_{n}$ is termed as
    $$
        \tau_{n} \equiv [\xi_n; \xi_{n+1}; \ldots].
    $$
    The distribution of $\tau_{n}$ is continuous if and only if the distribution of $\tau_{n+1}$ is continuous as well. For $\tau_{n}$, its cumulative distribution function is given by
    
    $$
    F_{\tau_n}(x)=P(\xi_n=\frac{1}{2})\cdot P\left(\frac{1}{\tau_{n+1}+\frac{1}{2}}<x\right)+P(\xi_n=1)\cdot P\left(\frac{1}{\tau_{n+1}+1}<x\right)=$$
        $$=p_{(\frac{1}{2})n}\left(1-F_{\tau_{n+1}}\left(\frac{1}{x}-\frac{1}{2}\right)\right)+p_{(1)n}\left(1-F_{\tau_{n+1}}\left(\frac{1}{x}-1\right)\right)=$$
    $$=1-p_{(\frac{1}{2})n}F_{\tau_{n+1}}\left(\frac{1}{x}-\frac{1}{2}\right)-p_{(1)n}F_{\tau_{n+1}}\left(\frac{1}{x}-1\right).$$

    Thus, if $F_{\tau_{n}}(x)$ contains both a singular and an absolutely continuous component, the same property must hold for $F_{\tau_{n+1}}(x)$. Without loss of generality, assume that $p_{(\frac{1}{2})(jL+i)} = q_{(\frac{1}{2})i}$ for all $j \in \mathbb{Z}_{+}$ and $i \in \{1, 2, \ldots, L\}$.

    By Theorem \ref{teo2}, there exists a set $B$ of positive Lebesgue measure that contains no $A_2$-binary numbers. For each $x = [\alpha_1; \alpha_2; \ldots] \in B$, condition \eqref{lev2} is satisfied, and there exists a finite value $F_{\xi}'(x) > 0$. Consequently, we have
    $$
   \lim\limits_{n\rightarrow\infty}\frac{p_{(\alpha_1)1}\cdot p_{(\alpha_2)2}\cdot\ldots\cdot p_{(\alpha_n)n}}{\lambda(\Delta_n(\alpha_1;\ldots;\alpha_n))}=\lim\limits_{n\rightarrow\infty}\frac{P(\xi\in\Delta_n(\alpha_1;\ldots;\alpha_n))}
   {\lambda(\Delta_n(\alpha_1;\ldots;\alpha_n))}=F_\xi^\prime(x)>0\Rightarrow
   $$
    \begin{equation}
    \label{a1}
       \Rightarrow \lim\limits_{n \rightarrow \infty} \sqrt[n]{p_{(\alpha_1)1} \cdot \ldots \cdot p_{(\alpha_n)n}} = e^{-2G}.
    \end{equation}
    
    On the probability space $([\frac{1}{2}, 1], \mathcal{B}([\frac{1}{2}, 1]), \eta(\cdot))$, we define random variables $(\psi_{n})$ in the following way:
    $$
        \psi_{n}(x) =
        \begin{cases}
            \ln(p_{(\frac{1}{2})n}), & \text{if } x \in T^{-n}(\Delta_{1}(\frac{1}{2})), \\
            \ln(p_{(1)n}), & \text{otherwise}.
        \end{cases}
    $$
    The invariance of $\eta(\cdot)$ under $T$ implies
    $$
        P(\psi_{n} = \ln(p_{(\frac{1}{2})n})) = \eta(T^{-n}(\Delta_{1}(\frac{1}{2})))
    =\frac{\ln(x+1)-\ln(x+2)}{\ln\frac{10}{9}}\Bigg|_{0,5}^{\frac{2}{3}} = \frac{\ln\frac{25}{24}}{\ln\frac{10}{9}}.
    $$
    and hence 
    $$P(\psi_{n} = \ln(p_{(1)n})) = 1 - \frac{\ln\frac{25}{24}}{\ln\frac{10}{9}}.$$ 
   Next, we consider random variables $\theta_n$ given by
    $$
    \theta_{n}=\frac{1}{L}\left(\psi_{(n-1)L+1}+\psi_{(n-1)L+2}+...+\psi_{(n-1)L+L}\right) \;\; \forall n \in \mathbb{N}.
    $$

    Since the measures \( \eta(\cdot) \) and \( \lambda(\cdot) \) are equivalent, there exists a set $ B^{*} $ such that $ \eta(B^{*}) > 0 $ and, for every $ x \in B^{*} $, condition \eqref{a1} is satisfied. Consequently, from \eqref{a1}, it follows that for every $ x \in B^{*}$
    $$
    \lim_{n \to +\infty} \frac{1}{n} \sum_{k=1}^{n} \psi_{k}(x) = -2G \quad \Rightarrow \quad \lim_{n \to +\infty} \frac{1}{n} \sum_{k=1}^{n} \theta_{k}(x) = -2G.
    $$
    The final condition, according to the Hewitt-Savage law, can occur only for events with probability 0 or 1. Since $ \eta(B^{*}) > 0 $, the condition  
$$
\lim_{n \to +\infty} \frac{1}{n} \sum_{k=1}^{n} \theta_{k} = -2G
$$
is satisfied with probability 1.

Taking into account the Lebesgue Dominated Convergence Theorem, we have  
\[
\sum_{j=1}^{L} \left( \frac{\ln\frac{25}{24}}{\ln\frac{10}{9}} \ln(q_{(0.5)j}) + \left(1 - \frac{\ln\frac{25}{24}}{\ln\frac{10}{9}}\right) \ln(q_{(1)j}) \right) = 
\lim_{n \to +\infty} M\left(\frac{1}{n} \sum_{k=1}^{n} \theta_{k}\right) = -2G,
\]
which contradicts condition \eqref{main}.

\textbf{Subcase 2:} $q_{(\frac{1}{2})j} q_{(1)j} = 0$ for some $j \in \{1, 2, \ldots, L\}$. Without loss of generality, assume $q_{(\frac{1}{2})1} = 0$. Let us consider a function
    $$
        f(x) =
        \begin{cases}
            1, & x \in \Delta_{2L}(\underbrace{\frac{1}{2}; \ldots; \frac{1}{2}}_{2L}), \\
            0, & x\in [\frac{1}{2}, 1]\setminus \Delta_{2L}(\underbrace {\frac{1}{2},\frac{1}{2},..., \frac{1}{2}}_{2L}).
        \end{cases}
    $$
    By the Birkhoff-Khinchin theorem, for almost all numbers $x \in [\frac{1}{2}, 1]$ in the sense of the measure $\eta(\cdot)$ (and consequently $\lambda(\cdot)$), we have:  
    $$
        \lim_{n \to \infty} \frac{f(T(x)) + \ldots + f(T^{n}(x))}{n} = \displaystyle\int\limits_{\Delta_{2L}(\underbrace{\frac{1}{2}; \ldots; \frac{1}{2}}_{2L})} \frac{dx}{\ln\frac{10}{9}(x+1)(x+2)}.
    $$
    
    It follows that for almost all $x \in [\frac{1}{2}, 1]$, the digit collection $(\frac{1}{2}, \frac{1}{2}, \ldots, \frac{1}{2})$ with $2L$ repetitions appears infinitely often in their $A_{2}$-representation. Consequently, condition $\eqref{a1}$ cannot hold on a positive Lebesgue measure set. This leads to a contradiction, and hence the distribution of $\xi$ is singular.

    Thus, in both subcases, the distribution of $\xi$ is singular.
\end{proof}

Let us summarize everything discussed so far: we have identified the necessary and sufficient conditions for the discreteness of the distribution of $\xi$, as well as partial conditions for its singularity. However, the type of distribution of $\xi$ remains unknown in more general settings. Another interesting question is how the distribution depends on the choice of the alphabet in the $A_2$-representation.

\begin{acknowledgement}
The first and second authors were supported by a grant from the Simon Foundation (1290607, P. M. V. and M. O. P.). The third author would like to thank the Isaac Newton Institute for Mathematical Sciences and the London Mathematical Society for their comprehensive support, as well as the University of St. Andrews for their hospitality.
\end{acknowledgement}

\ethics{Competing Interests}{The authors have no conflicts of interest to declare that are relevant to the content of this chapter.}

%
%
%


\end{document}